%% file: Vsyo_otnositel_no_.tex
\documentclass[11pt,reqno]{amsart}
\usepackage{amssymb,amscd,amsbsy}
\usepackage{amssymb,amscd,amsbsy,mathrsfs}
\input{classstartscr}

\theoremstyle{remark}

\newtheorem*{rem*}{Remark}

\newcommand{\OL}{{\rm OL}}

\newcommand{\MOL}{{\frak M}_{\OL}}

\newcommand\fM{\frak M}

\newcommand\dg{\frak D}

\newcommand{\rd}{{\rm d}}

\newcommand\ri{{\rm i}}
\newcommand\x{{\bs x}}
\newcommand\y{{\bs y}}
\newcommand\bx{{\bf x}}
\newcommand\rI{{\rm I}}
\newcommand\II{{\rm II}}

\begin{document}

\newcommand{\vse}{\vspace{.2in}}
\numberwithin{equation}{section}

%\title{Functions of self-adjoint operators under relatively bounded and relatively trace class perturbations. Relatively operator Lipschitz functions}
\author{A.B. Aleksandrov and V.V. Peller}

\title{Functions of self-adjoint operators under relatively bounded and relatively trace class perturbations}
%\author{V.V. Peller}
\thanks{The research on \S\:3-5  is supported by 
Russian Science Foundation [grant number 23-11-00153].
The research on the rest of the paper is supported by a grant of the Government of the Russian Federation for the state support of scientific research, carried out under the supervision of leading scientists, agreement  075-15-2024-631.}
%\thanks{Corresponding author: V.V. Peller; email: peller@math.msu.edu}

\begin{abstract}
We study the behaviour of functions of self-adjoint operators under relatively bounded and relatively trace class perturbation. We introduce and study the class of relatively operator Lipschitz functions. An essential role is played by double operator integrals. We also consider the class of resolvent Lipschitz functions. 
Then we obtain a trace formula in the case of relatively trace class perturbations and show that the maximal class of function for which the trace formula holds in the case of relatively trace class perturbations coincides with the class of relatively operator Lipschitz functions. Our methods also gives us a new approach to the inequality $\int|\bs{\xi}(t)|(1+|t|)^{-1}\,{\rm d}t<\be$ for the spectral shift function $\bs\xi$  in the case of relatively trace class perturbations.
\end{abstract}

\maketitle

{\bf
\footnotesize
\tableofcontents
\normalsize
}

\setcounter{section}{0}
\section{\bf Introduction}
\setcounter{equation}{0}
\label{In}

\

In this paper we are going to study the behaviour of functions of self-adjoint operators under relatively bounded and relatively trace class perturbations.

Let us start, however, with the classical case of bounded and trace class perturbations. Let $A$ be a (possibly unbounded) self-adjoint operator and let $K$ be a bounded self-adjoint operator. The operator $A$ can be interpreted as the initial (or unperturbed operator), $K$ is considered as a perturbation while the operator $B=A+K$ is the perturbed operator.

For a function $f$ continuous on the real line $\R$, the problem is to estimate the size (in one or another sense)
of the operator $f(B)-f(A)$ in terms of the size of $B-A$. A function $f$ is said to be an {\it operator Lipschitz function} if
\bay
\label{opLip}
\|f(B)-f(A)\|\le\const\|B-A\|
\ey
whenever $B-A$ is a bounded operator. Note that if inequality \rf{opLip} holds for all bounded self-adjoint operator $A$ and $B$, then it also holds for all not necessarily bounded operators.

The class of operator Lipschitz functions on $\R$ will be denoted by
$\OL(\R)$ or simply $\OL$ if it is clear that one deals with functions on $\R$.

Similarly, if instead of self-adjoint operators, we consider unitary operators, we can define the class of operator Lipschitz functions on the unit circle (see \cite{APol}).

We refer the reader to the survey article \cite{APol}, which contains a lot of information about operator Lipschitz functions. Let us mention here that operator Lipschitz functions must be differentiable everywhere on $\R$
and differentiable at infinity which means that the limit
$$
\lim_{|t|\to\be}\frac{f(t)}t
$$
exists. However, operator Lipschitz functions do not have to be continuously differentiable, see \cite{APol}.

It is well known (see \cite{APol}) that $f$ is operator Lipschitz if and only if the following implication holds
\bay
\label{sokhryad}
B-A\in\bS_1\quad\Longrightarrow\quad f(B)-f(A)\in\bS_1,
\ey
where $\bS_1$ is the trace class. We refer the reader to \cite{GK} and \cite{BSSst} for the definition and basic properties of Schatten--von Neumann classes $\bS_p$.
Moreover \rf{sokhryad} is equivalent to the fact that $f$ is {\it trace class Lipschitz}, i.e.,
$$
\|f(B)-f(A)\|_{\bS_1}\le\const\|B-A\|_{\bS_1},
$$
whenever $B-A\in\bS_1$.

In \cite{L} physicist I.M. Lifshits when studying problems in quantum statistics and crystals theory
arrived at the problem to evaluate the trace of the operator difference $f(B)-f(A)$ in the case when $B-A$ is a perturbation of $A$ that is in a sense small compared to $A$. Mostly he considered the case when $B-A$ has finite rank. He discovered that there exists a real integrable function $\bs\xi$ on $\R$ that is determined by $A$ and $B$ such that
\bay
\label{sledy}
\trace\big(f(B)-f(A)\big)=\int_\R f'(t)\bs\xi(t)\,\rd t
\ey
for sufficiently nice functions $f$. The function $\bs\xi$ is unique and is called the {\it spectral shift function} associated with
$A$ and $B$.

Later M.G. Krein has given in \cite{Kr} a mathematically rigorous justification of the results of Lifshits and established that in the case when $B-A$ is a trace class operator the spectral shift function $\bs\xi$ is integrable. In \cite{Kr} a problem was posed to describe the maximal class of functions $f$, for which the above Lifshits--Krein trace formula holds for arbitrary self-adjoint operators $A$ and $B$ with trace class difference. The problem was solved in \cite{PeKr}. It turned out that the maximal class of such functions coincides with the class of operator Lipschitz functions.

Another classical situation is when instead of measuring the size of the perturbation in terms of the differences
$B-A$ we measure the size of the perturbation in terms of the difference of the resolvents 
$(B+\ri I)^{-1}-(A+\ri I)^{-1}$. 

A function $f$ on $\R$ is called a {\it resolvent Lipschitz function} if
$$
\|f(B)-f(A)\|\le\const\|(B+\ri I)^{-1}-(A+\ri I)^{-1}\|,
$$
for arbitrary self-adjoint operators $A$ and $B$.

We say that a self-adjoint operator $K$ is a {\it resolvent trace class perturbation}
of a self-adjoint operator $A$ if 
$$
(A+K+\ri I)^{-1}-(A+\ri I)^{-1}\in\bS_1.
$$

It is also well known that a function $f$ on $\R$ is a resolvent Lipschitz function if and only if
$$
f(B)-f(A)\in\bS_1,
$$
whenever $B-A$ is a resolvent trace class perturbation of $A$. Moreover, this is equivalent to the inequality
$$
\|f(B)-f(A)\|_{\bS_1}\le\const\|(B+\ri I)^{-1}-(A+\ri I)^{-1}\|_{\bS_1}.
$$
In the case of resolvent trace class perturbation there is an analog of the Lifshits--Krein trace formula. 
Indeed, suppose that $K$ is a resolvent trace perturbation of $A$ and $B=A+K$. Then there exists a real measurable function $\bs\xi$ that satisfies
$$
\int_\R\frac{|\bs\xi(t)|}{1+t^2}\,\rd t<\be
$$
and such that trace formula \rf{sledy} holds for sufficiently nice functions $f$. Note that unlike the case of trace class perturbations a spectral shift function $\bs\xi$ in the case of resolvent trace class perturbations is not unique. It is determined by $A$ and $B$ modulo an additive real constant. 

The above facts are well known and can be obtained by applying the Cayley transform and reducing the situation to the case of unitary operator, see e.g., \cite{MNP}.

In this paper we mostly concentrate on the case of relatively bounded and relatively trace class perturbations;
these notions will be introduced in \S\;\ref{RBRTC}. We also introduce and study in \S\;\ref{RBRTC} the class of relatively operator Lipschitz functions, see \S\:\ref{RBRTC} for the definition

We obtain in \S\:\ref{SMdd} various descriptions of the class of relatively operator Lipschitz functions. We also find in \S\:\ref{SMdd} various new descriptions 
of the class of resolvent Lipschitz functions.

An important role will be played by double operator integrals. We give a brief introduction to double operator integrals in \S\;\ref{dvoinya} and introduce the notion of Schur multipliers with respect to spectral measures.
We also discuss in \S\;\ref{dvoinya} how one can express the operator difference $f(B)-f(A)$ in terms of double operator integrals for various classes of operators and functions.

In Section \ref{predstavili} we obtain a representation of the operator difference $f(B)-f(A)$ in the case of a relatively bounded perturbation in terms of a double operator integral.

Then in \S\;\ref{neobkhodimo} we show that a function $f$ on $\R$ is relatively operator Lipschitz if and only if
the function 
$$
(x,y)\mapsto\frac{f(x)-f(y)}{x-y}(x+\ri)
$$
on $\R^2$ is a Schur multiplier with respect to arbitrary spectral measures.

Next, in Section \S\:\ref{komutator?} we show that a function is relatively operator Lipchitz if and only if certain commutator estimates hold.

In \S\:\ref{differentsiruem} we obtain a formula for the derivative of the function
$$
t\mapsto f(A+t(B-A))
$$
in the strong operator topology. This formula will be used in \S\:\ref{fss} to obtain a trace formula in the case of relatively trace class perturbations. Our methods allows us to prove the inequality 
$$
\int_\R\frac{|\bs\xi(t)|}{1+|t|}\,\rd t<\be.
$$
for the spectral shift function $\bs\xi$ for the pair $\{A,B\}$, where $B-A$ is a relatively trace class perturbation of $A$. This inequality was obtained earlier in \cite{CS}.
However, our trace formula is valid for a considerably broader class of functions compared to the results  of
\cite{CS}. Moreover, we describe the maximal class of functions, for which the trace formula holds in the case of relatively trace class perturbations and show that this class coincides with the class of relatively operator Lipschitz functions.

\

\section{\bf Relatively bounded and relatively trace class perturbations}
\setcounter{equation}{0}
\label{RBRTC}

\

Let us proceed now to relatively bounded perturbations and relatively trace class perturbations.

Let $A$ be a not necessarily bounded self-adjoint operator. Suppose that $K$ is a closed operator such that
\bay
\label{vklyu}
{\rm Dom}(A)\subset{\rm Dom}(K),
\ey
where the notation ${\rm Dom}(T)$ stands for the domain of an operator $T$. It is well known (see \cite{BSSst}, Theorem 1, Ch. 3, Sec. 4) that in this case
$K$ must be {\it dominated by} $A$, i.e., the following inequality must hold:
\bay
\label{otogner}
\|Kv\|\le\const(\|v\|+\|Av\|),\quad v\in{\rm Dom}(A).
\ey

If, in addition to this, $K$ is a self-adjoint operator, we say that
$K$ is a {\it relatively  
bounded self-adjoint perturbation of} $A$. 

It is easy to see that property \rf{otogner} is equivalent to the fact that the operator
\bay
\label{C}
C\df K(A+\ri I)^{-1}
\ey
is bounded which, in turn, is equivalent to the fact that the operator
\bay
\label{G}
G\df K(A^2+I)^{-1/2}
\ey
is bounded.

\medskip

%{\bf Remark.} It is easy to see that if $K$ is a self-adjoint operator satisfying
%the operator $(B-A)(A+\ri I)^{-1}$ is bounded if and only if the operator

\medskip

We say that a self-adjoint operator $K$ 
is a {\it relatively trace class perturbation of} $A$ if the operator $C$ defined by \rf{C} belongs to $\bS_1$ which is, clearly, equivalent to the fact that the operator $G$ defined by \rf{G} belongs to $\bS_1$.

It is well known and easy to verify that if $K$ is a relatively trace class perturbation of $A$, then it is a resolvent trace class perturbation of $A$. The converse is not true which can easily be seen from Theorems \ref{x+i}
and \ref{resLip}.

\medskip

We are going to make three observations. All three are well known facts.

\medskip

\begin{lem}
\label{pervaya}
Suppose that $K$ is a self-adjoint operator satisfying the inequality
\bay
\label{KcdA}
\|Kv\|\le c\|v\|+d\|Av\|,\quad v\in{\rm Dom}(A),\quad\mbox{for some}\quad c>0\quad\mbox{and}\quad d\in(0,1)
\ey
(such operators are called {\it strictly dominated by $A$}). Then $K$ is also dominated by $A+tK$ for
$t\in(0,1]$.
\end{lem}

\Pf
Indeed, let $v\in{\rm Dom}(A)$. We have
$$
\|Kv\|\le c\|v\|+d\|Av\|=c\|v\|+d\|(A+tK)v-tKv\|\le c\|v\|+d\|(A+tK)v\|+dt\|Kv\|.
$$
It follows that
\bay
\label{cdK}
\|Kv\|\le(1-td)^{-1}\big(c\|v\|+d\|(A+tK)v\|\big),
\ey
and so $K$ is dominated by $A+tK$. $\bl$

\medskip

\begin{lem}
\label{vtoraya}
Suppose that $K$ is a self-adjoint operator that is a {\it relatively compact perturbation of} $A$, i.e., the operator $K(A+\ri I)^{-1}$ is compact. Then for any positive number $d$, there exists a positive number $c$ such that
$$
\|Kv\|\le c\|v\|+d\|Av\|,\quad v\in{\rm Dom}(A),
$$
and so, $K$ is strictly dominated by $A+K$.
\end{lem}

\Pf Let $\e>0$. By the assumption, the operator $K(A+\ri I)^{-1}$ is compact, and so it can be represented in the form $K(A+\ri I)^{-1}=T+R$, where $T$ is a finite rank operator and $\|R\|<\e$. 

Consider the operators $K_1$ and $K_2$ defined on ${\rm Dom}(A)$ by
$$
K_1=T(A+\ri I)\quad\mbox{and}\quad K_2=R(A+\ri I).
$$
Clearly,
$$
\|K_2v\|=\|R(A+\ri I)v\|\le\e\|(A+\ri I)v\|\le\e\|Av\|+\e\|v\|,\quad v\in{\rm Dom}(A),
$$
and so, it suffices to show that for any positive number $d$, there exists a positive number $c$ such that
$$
\|K_1v\|\le c\|v\|+d\|Av\|,\quad v\in{\rm Dom}(A).
$$
Since $T$ is a linear combination of rank one operators, it suffices to consider the case when 
$\rank T=\rank\big(K_1(A+\ri I)^{-1}\big)=1$. Suppose that
$$
K_1v=\big((A+\ri I)v,\f\big)\psi,\quad v\in{\rm Dom}(A).
$$
since $A$ is densely defined, $\f$ admit a representation $\f=\f_1+\f_2$, where 
$\f_1\in{\rm Dom}(A)$ and $\|\f_2\|<\e\|\psi\|^{-1}$. We have
\begin{align*}
K_1v&=\big((A+\ri I)v,\f_1\big)\psi+\big((A+\ri I)v,\f_2)\psi\big)\\[.2cm]
&=\big(v,(A+\ri I)\f_1)\psi+\big((A+\ri I)v,\f_2\big)\psi,
\quad v\in{\rm Dom}(A),
\end{align*}
and so
\begin{align*}
\|Kv\|&\le\|(A+\ri I)\f_1\|\cdot\|\psi\|\cdot\|v\|+\|\f_2\|\cdot\|\psi\|\cdot\|v\|+\|\f_2\|\cdot\|\psi\|\cdot\|Av\|\\[.2cm]
&\le\big(\|(A+\ri I)\f_1\|\cdot\|\psi\|+\|\f_2\|\cdot\|\psi\|\big)\|v\|+\e\|Av\|,
\quad v\in{\rm Dom}(A).\quad\bl
\end{align*}

\begin{lem}
\label{tret'ya}
Suppose that a self-adjoint operator $K$ is a relatively trace class perturbation of a self-adjoint operator $A$.
Then $K$ is a relatively trace class perturbation of $A+K$.
\end{lem}

\Pf It is easy to see that
$$
K(A+K+\ri I)^{-1}-K(A+\ri I)^{-1}=-K(A+K+\ri I)^{-1}K(A+\ri I)^{-1}\in\bS_1.
$$
Indeed, by the assumption, $K(A+\ri I)^{-1}\in\bS_1$. On the other hand, by Lemmata \ref{pervaya} and
\ref{vtoraya}, the operator $K(A+K+\ri I)^{-1}$ is bounded. $\bl$

\medskip

{\bf Definition.}
A bounded continuous function $f$ on $\R$ is called {\it relatively operator Lipschitz} if there is a positive number $k$ such that
\bay
\label{rOL}
\|f(B)-f(A)\|\le k\|(B-A)(A+\ri I)^{-1}\|
\ey
whenever $A$ and $B$ are self-adjoint operators such that $B-A$ is a relatively bounded perturbation of $A$.
{\it We denote by {\rm ROL} the class of relatively operator Lipschitz functions on} $\R$.

\medskip

Note that the assumption that $f$ is bounded is natural. Indeed, if $f$ is a continuous function on $\R$ ans
we apply inequality \rf{rOL} to self-adjoint operators on the one-dimensional Hilbert space, we obtain the inequality
\bay
\label{neobkhodimo}
|f(s)-f(t)|\le\const|s-t|\cdot(1+|t|)^{-1},\quad s,~t\in\R,
\ey
which implies that $f$ is a bounded function.

\medskip

{\bf Remark 1.} For a function $f$ on $\R$ to be relatively operator Lipschitz, it suffices to show that inequality
\rf{rOL} holds for {\it bounded} operators $A$ and $B$. This can be proved in the same way as the corresponding fact for operator Lipschitz functions, see \cite{APol}. Indeed, this follows from Theorem 3.2.1 of
\cite{APol}.

\medskip 

{\bf Remark 2.} Obviously, if $f$ is a relatively operator Lipschitz function, then $f$ is an operator Lipschitz function, and so, $f$ is differentiable on $\R$. Indeed, \rf{rOL} implies that
$$
\|f(B)-f(A)\|\le k\|(B-A)(A+\ri I)^{-1}\|\le k\|A-B\|\cdot\|(A+\ri I)^{-1}\|\le k\|A-B\|.
$$

\medskip

The following theorem gives an elementary necessary condition for a function on $\R$ to be relatively operator Lipchitz.

\begin{thm}
\label{nonedostatochno}
Let $f$ be a relatively operator Lipschitz function on $\R$. Then the following equivalent condition holds:

{\em(a)} 
$$
|f(s)-f(t)|\le\const|s-t|\cdot(1+|t|)^{-1},\quad s,~t\in\R;
$$

%{\em(a)}  there exists a constant $c$ such that
%$|f(x)-f(y)|\le c\frac{|x-y|}{1+|x|}$ for all $x,y\in\R$;
%
{\em(b)}
$$
|f(s)-f(t)|\le\const\frac{|s-t|}{1+|s|+|t|},\quad s,~t\in\R;
$$

{\em(c)}
the function
$$
s\mapsto(s+\ri)f(s),\quad s\in\R,
$$
is a Lipschitz function on $\R$.
\end{thm}

\Pf We have already mentioned that (a) is a necessary condition, see \rf{neobkhodimo}.  \rf{neobkhodimo} 

The implication (b)$\Longrightarrow$(a) is trivial. To prove the implication  (a)$\Longrightarrow$(b), we
observe that (a) implies that $|f(s)-f(t)|\le\const|s-t|(1+|t|)^{-1}$ for $s,~t\in\R$,
and \lb$\min((1+|s|)^{-1},(1+|t|)^{-1})\le2(1+|s|+|t|)^{-1}$.

It remains to prove that (a)$\Longleftrightarrow$(c).
Substituting $t=0$ we obtain that each of statements (a) and (c)
implies that $f$ is bounded.

First we prove that (a)$\Longrightarrow$(c). We have
\begin{align*}
|(s+\ri)f(s)-(t+\ri)f(t)|&=|(s+\ri)(f(s)-f(t))+f(t)(s-t)|\\[.2cm]
&\le
|(s+\ri)(f(s)-f(t))|+|f(t)(s-t)|\\[.2cm]
&\le\const\frac{|s+\ri|}{1+|s|}|s-t|+\sup_{u\in\R}|f(u)|\cdot|s-t|.
\end{align*}

Similarly, we can prove that (c)$\Longrightarrow$(a):
\begin{align*}
|(s+\ri)(f(s)-f(t))|&=|(s+\ri)f(s)-(t+\ri)f(t)-f(t)(s-t)|\\[.2cm]
&\le|(s+\ri)f(s)-(t+\ri)f(t)|+|f(t)(s-t)|\\[.2cm]
&\le\const|s-t|+\sup_{u\in\R}|f(u)|\cdot|s-t|.
\end{align*}
It remains to observe that $|s+\ri|^{-1}\le\frac{\sqrt 2}{1+|s|}$. $\bl$

\medskip

In \S\;\ref{neobkhodimost'} we observe that condition \rf{neobkhodimo} is not sufficient for a function to be relatively operator Lipsschitz.

The following fact can proved in the same way as the corresponding result for operator Lipschitz function, see
Theorem 3.6.5 in \cite{APol}.

\begin{thm}
\label{rolryal}
Let $f$ be a function on $\R$. The following statements are equivalent:

{\rm(a)} $f\in{\rm ROL}$; 

{\rm(b)} if $K$ is a self-adjoint relatively trace class perturbation of a self-adjoint operator $A$, then
$f(A+K)-f(A)\in\bS_1$;

{\rm(c)} $f$ is a relatively trace class Lipschitz function, i.e.,
$$
\|f(A+K)-f(A)\|_{\bS_1}\le\const\|K(A+\ri I)^{-1}\|_{\bS_1}
$$
for an arbitrary relatively trace class perturbation $K$ of a self-adjoint operator $A$.
\end{thm}

\

\section{\bf The role of double operator integrals and Schur multipliers}
\setcounter{equation}{0}
\label{dvoinya}

\

Double operator integrals are expressions of the form
\bay
\label{dvopi}
\iint\Phi(x,y)\,\rd E_1(x)Q\,\rd E_2(y).
\ey
Here $E_1$ and $E_2$ are spectral measures on Hilbert space, $\Phi$ is a bounded measurable function and 
$Q$ is a bounded linear operator on Hilbert space.

Double operator integrals play a very important role in perturbation theory. They appeared first in the paper by Yu.L. Daletskii and S.G. Krein \cite{DK}. Later M.S. Birman and M.Z. Solomyak created in \cite{BS1,BS2}
and \cite{BS4} a beautiful rigorous theory of double operator integrals.

The starting point of the Birman--Solomyak theory is the case when $Q$ belongs to the class $\bS_2$ of Hilbert--Schmidt operators. Consider the set function $\E$ that takes values in the set of orthogonal projections on the Hilbert Schmidt class $\bS_2$ and is defined on the measurable rectangles by
$$
\E(\L\times\D)T\df E_1(\L)TE_2(\D),\quad T\in\bS_2.
$$
Birman and Solomyak showed in \cite{BStp} that $\E$ extends to a spectral measure on $\bS_2$. This allows one to define the double integral \rf{dvopi} in the case when $Q\in\bS_2$ by
$$
\iint\Phi\,\rd E_1Q\,\rd E_2\df\left(\int\Phi\,\rd\E\right)Q.
$$

Next, $\Phi$ is said to be a {\it Schur multiplier with respect to $E_1$ and} $E_2$ if
$$
Q\in\bS_1\quad\Longrightarrow\quad\iint\Phi\,\rd E_1Q\,\rd E_2\in\bS_1.
$$
We use the notation $\fM_{E_1,E_2}$ for the class of Schur multipliers with respect to $E_1$ and $E_2$.
The norm $\|\Phi\|_{\fM_{E_1,E_2}}$ in the space of Schur multipliers is defined as the norm of the transformer
\bay
\label{transformator}
Q\mapsto\iint\Phi\,\rd E_1Q\,\rd E_2
\ey
on $\bS_1$.

In the case when $\Phi\in\fM_{E_1,E_2}$, one can define by duality the double operator integral \rf{dvopi}
for arbitrary bounded operator $Q$. Moreover, the norm of the transformer \rf{transformator} on the space of bounded linear operators also coincides with $\|\Phi\|_{\fM_{E_1,E_2}}$.

There are several characterizations of the class $\fM_{E_1,E_2}$, see, e.g., \cite{Pe1} and \cite{AP3}. In particular, $\Phi\in\fM_{E_1,E_2}$ if and only if $\Phi$ belongs to the {\it Haagerup tensor product}
$L^\be_{E_1}\!\otimes_{\rm h}L^\be_{E_2}$ of $L^\be$ spaces $L^\be_{E_1}$ and $L^\be_{E_2}$, i.e.,
$\Phi$ admits a representation
\bay
\label{predstavim}
\Phi(x,y)=\sum_n\f_n(x)\psi_n(y),
\ey
where $\f_n\in L^\be_{E_1}$, $\psi_n\in L^\be_{E_2}$, and
$$
\{\f_n\}_{n\ge0}\in L_{E_1}^\be(\ell^2)\quad\mbox{and}\quad
\{\psi_n\}_{n\ge0}\in L_{E_2}^\be(\ell^2).
$$
By the {\it norm of $\Phi$ in $L^\be(E_1)\!\otimes_{\rm h}\!L^\be(E_2)$} we mean the infimum of
\bay
\label{normaHaag}
\big\|\{\f_n\}_{n\ge0}\big\|_{L_{E_1}^\be(\ell^2)}
\big\|\{\psi_n\}_{n\ge0}\big\|_{L_{E_2}^\be(\ell^2)}
\ey
over all representations in the form \rf{predstavim}. Moreover, it was established recently in \cite{AP3}
that this description is {\it isometric}, i.e., the norm of $\Phi$ in $\fM_{E_1,E_2}$ coincides with the 
norm of $\Phi$ in $L^\be_{E_1}\!\otimes_{\rm h}L^\be_{E_2}$.

In the case when $\Phi\in L^\be_{E_1}\!\otimes_{\rm h}L^\be_{E_2}$, the double operator integral \rf{dvopi} can be computed as follows:
$$
\iint\Phi\,\rd E_1Q\,\rd E_2=\sum_n\left(\int\f_n\,\rd E_1\right)Q\left(\int\psi_n\,\rd E_2\right)
$$
and the integral on the right converges in the weak operator topology.

We refer the reader to the survey article \cite{APol} for more detailed information on double operator integrals and Schur multipliers.

It turns out that a function $f$ on $\R$ is operator Lipschitz if and only if it is differentiable on $\R$ and the divided difference $\dg f$ defined by
\bay
\label{razdraz}
(\dg f)(x,y)\df\left\{\begin{array}{ll}\dfrac{f(x)-f(y)}{x-y},&\text {if}\,\,\,\,x\ne y,\\[.4cm]
f'(x),&\text{if}\,\,\,\,x=y,
\end{array}\right.
\ey
is a Schur multiplier with respect to arbitrary Borel spectral measures.\footnote{We use the notation
$\fM(\R^2)$ for the class of functions on $\R^2$ that are Schur multiplier with respect to arbitrary Borel spectral
measures.} Moreover, for operator Lipschitz functions $f$ the following representation holds for arbitrary self-adjoint operator $A$ and $B$ with bounded 
$B-A$:
\bay
\label{razndvoi}
f(B)-f(A)=\iint_{\R\times\R}(\dg f)(x,y)\,\rd E_B(x)(B-A)\,\rd E_A(y),
\ey
where $E_A$ and $E_B$ are the spectral measures of $A$ and $B$.

We refer the reader to \cite{APol} for more detailed information.

The class of resolvent Lipschitz functions also admits a characterization in terms of Schur multipliers: a differentiable function $f$ on $\R$ is a resolvent Lipschitz function if and only if the function
$$
(x,y)\mapsto(\dg f)(x,y)(x+\ri)(y+\ri)
$$
on $\R^2$ is a Schur multiplier with respect to arbitrary Borel spectral measures. The following formula
holds
$$
f(B)-f(A)=\iint(\dg f)(x,y)(x+\ri)(y+\ri)\,\rd E_B(x)\big((A+\ri I)^{-1}-(B+\ri I)^{-1}\big)\rd E_A(y)
$$
for resolvent Lipchitz functions $f$.
This can easily be reduced to an analogue of formula \rf{razndvoi} for functions of unitary operators by passing to Cayley transform.

It is also well known (see e.g., \cite{MNP}) that a function $f$ on $\R$ is resolvent Lipschitz if and only if
is the function $\f$ on $\T\setminus\{1\}$ defined by
\bay
\label{na okruzhnosst'}
\f(\z)=f\left(\ri\frac{1+\z}{1-\z}\right),
\ey
then $\f$ extends to an operator Lipschitz function on $\T$. This can also be proved by passing to Cayley transform.

Suppose now that $A$ and $B$ are self-adjoint operators and $B-A$ is a resolvent trace class perturbation of $A$. Then the trace formula holds
$$
\trace\big(f(B)-f(A)\big)=\int_\R f'(t)\bs\xi(t)
$$
for an arbitrary resolvent Lipschitz function $f$, where $\bs\xi$ is a {\it spectral shift function} that corresponds to $A$ and $B$. It is real valued, unique modulo a constant additive and satisfies the condition
$$
\int_\R\frac{|\bs\xi(t)|}{1+t^2}\,\rd t<\be.
$$
Moreover, the class of resolvent Lipschitz functions is the maximal class of functions, for which trace formula
\rf{razndvoi} holds as soon as $B-A$ is a resolvent trace class perturbation of $A$.

This again can be obtained by passing to the Cayley transforms of $A$ and $B$. We refer the reader to \cite{MNP} for details.

In this paper we characterize the class of relatively operator Lipschitz functions as the class of function $f$ on 
$\R$, for which the function
$$
(x,y)\mapsto(\dg f)(x,y)(x+\ri)
$$
on $\R^2$ is a Schur multiplier with respect to arbitrary Borel spectral measures.
In the case when $B-A$ is a relatively bounded perturbation of $A$ and $f$ is relatively operator Lipschitz, we  represent $f(B)-f(A)$ in terms of the double operator integral
$$
f(B)-f(A)=\iint(\dg f)(x,y)(y+\ri)\,\rd E_B(x)(B-A)(A+\ri I)^{-1}\rd E_A(y).
$$

\

\section{\bf Schur multipliers originating from divided differences}
\setcounter{equation}{0}
\label{SMdd}

\

Let $f$ be a differentiable function on $\R$. As before, we use the symbol $\dg f$
the divided difference \rf{razdraz}.

%Put
%$$
%(\dg_0 f)(x,y)\df\left\{\begin{array}{ll}\dfrac{f(x)-f(y)}{x-y},&\text {if}\,\,\,\,x\ne y,\\[.2cm]
%0,&\text{if}\,\,\,\,x=y.
%\end{array}\right.
%$$

\medskip

{\bf Definition.} We say that a function $\f$ on $\R$ is a multiplier of the space $\OL$ of operator Lipschitz functions if
$$
f\in\OL\quad\Longrightarrow\quad \f f\in\OL.
$$
We denote by $\MOL$ the class of all multipliers of $\OL$.
By the norm $\|\f\|_{\MOL}$ of a function $\f$ in $\MOL$ we mean the norm of the multiplication operator 
$f\mapsto\f f$ on the space $\OL$.

\medskip

For convenience, throughout the paper we introduce the following notation. We denote by $\x$ and $\y$ the functions on $\R^2$ defined by
$$
\x(s,t)\df s,\quad(s,t)\in\R^2,\quad\mbox{and}\quad \y(s,t)\df t,\quad(s,t)\in\R^2.
$$
We are also going to use the notation $\bx$ for the function on $\R$ defined by
$$
\bx(s)=s,\quad s\in\R.
$$

We start with analyzing the condition
$$
(\x+\ri)({\dg f})\in\fM(\R^2).
$$
As we have mentioned in \S\:\ref{dvoinya}, we are going to prove in this paper that this condition is equivalent to the fact that $f\in{\rm ROL}$, i.e., $f$ is a relatively operator Lipschitz function.

\begin{thm}
\label{x+i}
Let $f$ be a differentiable function on $\R$. The following are equivalent:

{\em(a)} 
$(\x+\ri)({\dg f})\in\fM(\R^2)$; 

{\em(b)}
$(\bx+\ri)f\in\OL(\R)$;

{\em(c)} $f\in\MOL$.
\end{thm}

\medskip

\Pf Let us first show that (a)$\Rightarrow$(b).
Suppose that $(\x+\ri)(\dg f)\in\fM(\R^2)$. Analizing $((\x+\ri)(\dg f))(x,y)$ for
$y=0$, we find that $f$ is bounded.
It remains to observe that
\bay
\label{ix}
\frac{(x+\ri)f(x)-(y+\ri)f(y)}{x-y}=(x+\ri)(\dg f)(x,y)+f(y),\quad x,~y\in\R.
\ey
Let us now prove that
(b)$\Rightarrow$(a).
Since $(\bx+\ri)f\in\OL(\R)$, it follows that 
$$
|(t+\ri)f(t)-if(0)|\le C|t|\quad\mbox{for every}\quad t\in\R.
$$
Consequently, 
$f$ is bounded. It follows now from \rf{ix} that $(\x+\ri)(\dg f)\in\frak M(\R^2)$. 

Next, the implication (c)$\Rightarrow$(b) is trivial because $(\bx+\ri)\in\OL(\R)$. It remains to establish that
(b)$\Rightarrow$(c).

Clearly, we may assume that $f(0)=0$.
Since $(\bx+\ri)f\in\OL(\R)$, it is easy to see that 
$$
|(t+\ri)f(t)|\le\|(\bx+\ri)f\|_{\OL(\R)}|t|\quad\mbox{for every}\quad t\in\R.
$$ 
Consequently, 
\bay
\label{fx+i}
|f(t)|<\|(\bx+\ri)f\|_{\OL(\R)}\quad\mbox{for every}\quad t\in\R.
\ey
Let us show that $fh\in\OL(\R)$ for an arbitrary function $h$ in $\OL(\R)$.

Let $A$ and $B$ be bounded self-adjoint operators.

Put $g=(\bx+\ri)^{-1}h$.
The identity
\begin{align*}
f(A)h(A)-f(B)h(B)&=f(A)(h(A)-h(B))+(f(A)-f(B))h(B)\\[.2cm]
&=f(A)(h(A)-h(B))+(f(A)(B+\ri I)-f(B)(B+\ri I))g(B)\\[.2cm]
&=f(A)(h(A)-h(B))+(f(A)(A+\ri I)-f(B)(B+\ri I))g(B)\\[.2cm]
&-f(A)(A-B)g(B)
\end{align*}
implies that
\begin{align*}
\|f(A)h(A)&-f(B)h(B)\|\le\|f\|_{L^\infty(\R)}\|h\|_{\OL(\R)}\|A-B\|\\[.2cm]
&+\|(\bx+\ri)f\|_{\OL(\R)}\|g\|_{L^\infty(\R)}\|A-B\|+\|f\|_{L^\infty(\R)}\|g\|_{L^\infty(\R)}\|A-B\|.
\end{align*}
By \rf{fx+i}, $\|f\|_{L^\infty(\R)}\le\|(\bx+\ri)f\|_{\OL(\R)}$. It remains to observe that $g\in L^\infty(\R)$.
Indeed,
$$
|g(t)|\le|t+\ri |^{-1}(|h(t)-h(0)|+|h(0)|)\le\|h\|_{\OL(\R)}+|h(0)|,\quad t\in\R. \quad\bl
$$

\medskip

Note that the equivalence of (b) and (c) in Theorem \ref{x+i} is a special case of Theorem 4.7  in \cite{A}.

\medskip

{\bf Remark.} Obviously, condition (a) in the statement of the theorem is equivalent to the fact that
$(\y+\ri)(\dg f)\in\fM(\R^2)$.

\medskip

Let us proceed now to the condition 
$$
(\x+\ri)(\y+\ri)(\dg f)\in\frak M(\R^2).
$$

As we have already mentioned in \S\:\ref{dvoinya}, this condition describes the class of resolvent Lipschitz functions. We have also mentioned in \S\:\ref{dvoinya} that this condition is equivalent to the fact that the function $\f$ on $\T\setminus\{1\}$ defined by \rf{na okruzhnosst'}
extends to an operator Lipschitz function on $\T$.

\begin{thm}
\label{resLip} 
Let $f$ be a differentiable function on $\R$. The following are equivalent:

{\em(a)}
$(\x+\ri)(\y+\ri)(\dg f)\in\frak M(\R^2)$;

{\em(b)}
$(\bx+\ri)^2(f-c)\in\OL(\R)$ for some $c\in\C$;

{\em(c)}
$(\bx+\ri)(f-c)\in\MOL$ for some $c\in\C$.

%{\em(d)} the function $\f$ on $\T\setminus\{1\}$ defined by \rf{na okruzhnosst'}
%extends to an operator Lipschitz function on $\T$.

If $c\in\C$ satisfies {\em(b)} or {\em(c)}, then $c$ must be equal to $\lim_{|x|\to\be}f(x)$.

\end{thm}

\Pf We need the following identity:
\begin{align}
\label{ixiy}
\frac{(x+\ri)^2f(x)-(y+\ri)^2f(y)}{x-y}&=(x+\ri)(y+\ri)(\dg f)(x,y)\nonumber\\[.2cm]
&+(x+\ri)f(x)+(y+\ri)f(y).
\end{align}
 
Let us show that (a)$\Rightarrow$(b).
Suppose that $(\x+\ri)(\y+\ri)(\dg f)\in\frak M(\R^2)$. Then \lb$(\x+\ri)(\dg f)\in\frak M(\R^2)$.
It follows from the previous theorem that $(\bx+\ri)f\in\OL(\R)$. Thus,
the limit $c\df\lim_{|x|\to\be}f(x)$ exists. Without loss of generality we may assume that $c=0$.
To deduce from \rf{ixiy} the inclusion $(\bx+\ri)^2f\in\OL(\R)$, it suffices to observe that
the function $(\bx+\ri)f$ is bounded. This follows from the fact that $\lim_{|x|\to+\infty}f(x)=0$ and from the fact 
that the function $(\bx+\ri)^2(\dg f)(\bx,\bx)=(\bx+\ri)^2f'$ is bounded.

Suppose now that $(\bx+\ri)^2(f-c)\in\OL(\R)$ for some $c\in\C$. Let us show that \lb$(\x+\ri)(\y+\ri)(\dg f)\in\frak M(\R^2)$. Again, we may assume that $c=0$. 
The result follows from the fact that the function $(\bx+\ri)f$ is bounded and from identity \rf{ixiy}. 

The equivalence of (b) and (c) is a consequence of in the statement of Theorem \ref{x+i}.

Finally, it is obvious that (b) can hold only for one value of $c$. The same is true for (c).
$\bl$

\medskip

{\bf Remark.} As we have already observed, condition (a) in the statement of Theorem \ref{resLip} is equivalent to the fact that
the function $\f$ on $\T\setminus\{1\}$ defined by \rf{na okruzhnosst'}
extends to an operator Lipschitz function on $\T$. This allows us to reduce the equivalence of
statements (a) and (b) to Theorem 5.6 of the paper \cite{A}.

%\begin{lem} 
%A function $w$ is a Schur multiplier of $\OL(\R)$  if and only if
% $(x+\ri)w(x)\in\OL(\R)$.
%\end{lem}
%
%\Pf The necessity of the condition $(x+\ri)w(x)\in\OL(\R)$ is obvious. Let us establish sufficiency.
%
%Clearly, we may assume that $w(0)=0$.
%Suppose that $(x+\ri)w(x)\in\OL(\R)$. Then
%$$
%|(t+\ri)w(t)|\le\|(x+\ri)w(x)\|_{\OL(\R)}|t|
%$$
%for every $t\in\R$.  
%Consequently, $|w(t)|<\|(x+\ri)w(x)\|_{\OL(\R)}$ for every $t\in\R$.
%Let us show that $wf\in\OL(\R)$ for an arbitrary function $f$ in $\OL(\R)$.
%
%Let $A$ and $B$ be bounded self-adjoint operators.
%
%Put $g(x)=(x+\ri)^{-1}f(x)$.
%The identity
%\begin{align*}
%w(A)f(A)-w(B)f(B)&=w(A)(f(A)-f(B))+(w(A)-w(B))f(B)\\[.2cm]
%&=w(A)(f(A)-f(B))+(w(A)(B+\ri I)-w(B)(B+\ri I))g(B)\\[.2cm]
%&=w(A)(f(A)-f(B))+(w(A)(A+\ri I)-w(B)(B+iI))g(B)\\[.2cm]
%&-w(A)(A-B)g(B)
%\end{align*}
%implies that
%\begin{align*}
%\|w(A)f(A)&-w(B)f(B)\|\le\|w\|_{L^\infty(\R)}\|f\|_{\OL(\R)}\|A-B\|\\[.2cm]
%&+\|(x+\ri)w(x)\|_{\OL(\R)}\|g\|_{L^\infty(\R)}\|A-B\|+\|w\|_{L^\infty(\R)}\|g\|_{L^\infty(\R)}\|A-B\|.
%\end{align*}
%It remains to observe that $\|w\|_{L^\infty(\R)}\le\|(x+\ri)w(x)\|_{\OL(\R)}$ and $g\in L^\infty(\R)$.
%Indeed,
%$$
%|g(x)|\le|x+\ri |^{-1}(|f(x)-f(0)|+|f(0)|)\le\|f\|_{\OL(\R)}+|f(0)|. \quad\bl
%$$

%\begin{thm} 
%The function
%$(x+\ri)(\dg f)(x,y)$ belongs to $\frak M(\R^2)$ if and only if $f$
%is a multiplier of $\OL(\R)$.
%\end{thm}
%
%\begin{thm} 
%The function
%$(x+\ri)(y+\ri)(\dg f)(x,y)$ belongs to $\frak M(\R^2)$ if and only if the function $(x+\ri)(f(x)-c)$
%is a multiplier of $\OL(\R)$
%for some $c\in\C$.
%\end{thm}

\

\section{\bf A representation in terms of double operator integrals of the operator\\ difference $\bs{f(B)-f(A)}$ in the case of relatively bounded perturbations}
\setcounter{equation}{0}
\label{predstavili}

\

\medskip

Suppose that $f$ is a differentiable function on $\R$. Consider the functions $\dg_\rI f$ and $\dg_\II f$ on $\R^2$ defined by
\bay
(\dg_\rI f)(x,y)\df\frac{f(x)-f(y)}{x-y}(y+\ri)\quad\mbox{and}\quad(\dg_\II f)(x,y)\df\frac{f(x)-f(y)}{x-y}(y^2+1)^{1/2}.
\ey
By the value of $(f(x)-f(y))(x-y)^{-1}$ in the case when $x=y$ we mean the derivative $f'(x)$.

It is easy to see that $\dg_\II f$ is a Schur multiplier with respect to arbitrary Borel spectral measures on $\R$ if and only if $\dg_\rI f$ has this property. 

The following fact was established in \cite{BS4}, Th. 4.6. We give a proof in this paper for completeness.

\begin{thm}
\label{otopLi}
Suppose that $A$ and $B$ are self-adjoint operators such that $B-A$ is a relatively bounded perturbation of $A$. Let $f$ be a differentiable function on $\R$ such that $\dg_\II f$ is a Schur multiplier with respect to arbitrary spectral measures on $\R$.
Then $f$ is a relatively operator Lipschitz function and
\begin{align}
\label{2ra}
f(B)&-f(A)=
\iint_{\R\times\R}\frac{f(x)-f(y)}{x-y}(y+\ri)\,{\rm d}E_{B}(x)(B-A)(A+\ri I)^{-1}\,{\rm d}E_A(y)\nonumber\\[.2cm]
&=\iint_{\R\times\R}\frac{f(x)-f(y)}{x-y}(y^2+1)^{1/2}\,{\rm d}E_{B}(x)(B-A)(A^2+I)^{-1/2}\,{\rm d}E_A(y).
\end{align}
\end{thm}

\Pf Clearly, the double operator integral
$$
\iint_{\R\times\R}\frac{f(x)-f(y)}{x-y}(y^2+1)^{1/2}\,{\rm d}E_{B}(x)(B-A)(A^2+I)^{-1/2}\,{\rm d}E_A(y)
$$
represents a bounded linear operator and 
\begin{multline}
\label{rOLDOI}
\left\|
\iint_{\R\times\R}\frac{f(x)-f(y)}{x-y}(y^2+1)^{1/2}\,{\rm d}E_{B}(x)(B-A)(A^2+I)^{-1/2}\,{\rm d}E_A(y)\right\|\\[.2cm]
\le\|\dg_\II\|_{\fM(\R^2)}\|(B-A)(A^2+I)^{-1/2}\|.
\end{multline}

For a positive number $M$, we consider the spectral projections, $P_M\df E_A[-M,M]$ and $Q_M\df E_{B}[-M,M]$. We are going to use the following agreement: by $\int_{-M}^M$ we mean $\int_{[-M,M]}$, i.e., the integral over the interval $[-M,M]$.

We have
\begin{multline*}
Q_M\left(
\iint_{\R\times\R}\frac{f(x)-f(y)}{x-y}(y^2+1)^{1/2}\,{\rm d}E_{B}(x)G(B-A)(A^2+I)^{-1/2}\,{\rm d}E_A(y)\right)P_M\\[.2cm]
=\int\limits_{-M}^M\int\limits_{-M}^M\frac{f(x)-f(y)}{x-y}(y^2+1)^{1/2}\,{\rm d}E_{B}(x)(B-A)(A^2+I)^{-1/2}\,{\rm d}E_A(y)\\[.2cm]
=\int\limits_{-M}^M\int\limits_{-M}^M\frac{f(x)-f(y)}{x-y}(y^2+1)^{1/2}\,{\rm d}E_{B}(x)(B-A)(A^2+I)^{-1/2}\,{\rm d}E_A(y)\\[.2cm]
=\int_{-M}^M\int_{-M}^M\frac{f(x)-f(y)}{x-y}(y^2+1)^{1/2}(y^2+1)^{-1/2}\,{\rm d}E_{B}(x)(B-A)\,{\rm d}E_A(y)
\\[.2cm]
=\int_{-M}^M\int_{-M}^M\frac{f(x)-f(y)}{x-y}\,{\rm d}E_{B}(x)(B-A)\,{\rm d}E_A(y).
\end{multline*}
Note that the operator $B$ is bounded on $\Range Q_M$ while $A$ is bounded on $\Range P_M$. Thus,
\begin{multline*}
\int\limits_{-M}^M\int\limits_{-M}^M\frac{f(x)-f(y)}{x-y}\,{\rm d}E_{B}(x)(B-A)\,{\rm d}E_A(y)=\\[.2cm]
=\int\limits_{-M}^M\int\limits_{-M}^M\frac{f(x)-f(y)}{x-y}\,{\rm d}E_{B}(x)B\,{\rm d}E_A(y)
-\int\limits_{-M}^M\int\limits_{-M}^M\frac{f(x)-f(y)}{x-y}\,{\rm d}E_{B}(x)A\,{\rm d}E_A(y)\\[.2cm]
=\int\limits_{-M}^M\int\limits_{-M}^Mx\frac{f(x)-f(y)}{x-y}\,{\rm d}E_{B}(x)\,{\rm d}E_A(y)-
\int\limits_{-M}^M\int\limits_{-M}^My\frac{f(x)-f(y)}{x-y}\,{\rm d}E_{B}(x)\,{\rm d}E_A(y)\\[.2cm]
=\!\!\int\limits_{-M}^M\int\limits_{-M}^M(x-y)\frac{f(x)-f(y)}{x-y}\,{\rm d}E_{B}(x)\,{\rm d}E_A(y)\!=\!
\int\limits_{-M}^M\int\limits_{-M}^M(f(x)-f(y))\,{\rm d}E_{B}(x)\,{\rm d}E_A(y)\\[.2cm]
=\int\limits_{-M}^M\int\limits_{-M}^Mf(x)\,{\rm d}E_{B}(x)\,{\rm d}E_A(y)
-\int\limits_{-M}^M\int\limits_{-M}^Mf(y)\,{\rm d}E_{B}(x)\,{\rm d}E_A(y)\\[.2cm]
=Q_M(f(B)-f(A))P_M.
\end{multline*}
Since both $Q_M$ and $P_M$ converge as $M\to\be$ to I in the strong operator topology, it follows that
$f(B)-f(A)$ is a bounded operator and $f(B)-f(A)$ is equal to the very last term in \rf{2ra}. The first equality in \rf{2ra} can be proved in exactly the same way. 

Finally, inequality \rf{rOLDOI} exactly means that $f$ is relatively operator Lipschitz.
$\bl$

\medskip

Below we prove Theorem \ref{Dfy+imS}, which asserts that the converse is also true, i.e., if $f$ is relatively operator Lipschitz, then $\dg_\rI f$ and $\dg_\II f$ belong to $\fM(\R^2)$.

\begin{cor}
Suppose that $K$ is a self-adjoint relatively trace class perturbation of a self-adjoint operator $A$ and let $f$ be a differentiable function such that $\dg_\rI f$ is a Schur multiplier with respect to arbitrary Borel spectral measures on 
$\R$. Then $f(A+K)-f(A)\in\bS_1$.
\end{cor}

\Pf Indeed, under the hypotheses of the corollary the transformer
$$
D\mapsto\iint_{\R\times\R}\frac{f(x)-f(y)}{x-y}(y+\ri)\,{\rm d}E_{A+K}(x)D\,{\rm d}E_A(y)
$$
is a bounded operator on $\bS_1$. The result follows from \rf{2ra}. $\bl$

Recall that it was established in \cite{Pe1} and \cite{Pe2} that functions in the Besov class $B_{\be,1}^1(\R)$ on the real line must be operator Lipschitz. This together with
Theorem \ref{otopLi} and Theorem \ref{x+i} allows us to obtain the following sufficient condition for relative operator Lipschitzness.

\begin{thm}
Let $f$ be a function on $\R$ such that the function $g$ defined by
$$
g(x)=f(x)(x+\ri),\quad x\in\R,
$$
belongs to the Besov class $B_{\be,1}^1(\R)$. Then $f$ is a relatively operator Lipschitz functions.
\end{thm}`

%\begin{thm}
%\label{obratno}
%Let $f$ be a relatively operator Lipschitz function on $\R$. Then the functions $\dg_\rI f$ and $\dg_\II f$ belong to 
%$\fM$, and so conditions {\rm(a)-(c)} in the statement of Theorem {\em\ref{x+i}} hold. 
%\end{thm}
%
%\Pf Let us first observe that $f$ is a bounded function. Indeed, we can apply inequality \rf{rOL} to operators on the one-dimensional Hilbert space. Without loss of generality, we may assume that $f(0)=0$. Then inequality \rf{rOL} implies that
%$$
%|f(t)|\le k\frac{|t|}{|t+\ri|},\quad t\in\R,
%$$
%and so $f$ is bounded.
%
%By Theorem \ref{x+i}, to prove the theorem, we can show that $(\bx+\ri)f$ is an operator Lipschitz function.
%Suppose that $A$ and $B$ are bounded self-adjoint operators. 
%We have
%$$
%f(B)(B+\ri I)-f(A)(A+\ri I)=f(B)(B-A)+\big(f(B)-f(A)\big)(A+\ri I).
%$$
%Clearly,
%$$
%\|f(B)(B-A)\|\le\|f\|_{L^\be(\R)}\|B-A\|.
%$$
%On the other hand, let $u$ be a vector in our Hilbert space. Put
%$v\df(A+\ri I)u$.
%$$
%\|\big(f(B)-f(A)\big)(A+\ri I)u\|=\|\big(f(B)-f(A)\big)v\|
%\le k\|(B-A)(A+\ri I)^{-1}\|\cdot\|v\|,
%$$
%where $k$ is the constant in \rf{rOL}. Thus,
%$$
%\|\big(f(B)-f(A)\big)(A+\ri I)u\|\le k\|(B-A)(A+\ri I)^{-1}\|\cdot\|(A+\ri I)u\|.
%$$
%%Since $A+\ri I$ is an invertible operator, the last inequality means that
%In other words,
%%$$
%%\|\big(f(B)-f(A)\big)\|\le k\|(B-A)(A+\ri I)^{-1}\|.
%%$$
%$$
%\|\big(f(B)-f(A)\big)v\|\le k\|(B-A)(A+\ri I)^{-1}\|\cdot\|v\|.
%$$
%It remains to observe that $\|(A+\ri I)^{-1}\|\le 1$. $\bl$

\

\section{\bf The necessity of the condition $\bs{\dg_\rI f\in\fM(\R^2)}$}
\setcounter{equation}{0}
\label{neobkhodimost'}

\

To establish the necessity of the condition $\dg_\rI f\in\fM(\R^2)$ for $f$ to be relatively operator Lipschitz,
we consider the one-parametric family $A_t\df A+tK$, $0\le t\le1$, and study the behaviour of the operators $f(A_t)-f(A)$. Here $A$ and $B$ are a self-adjoint operators on Hilbert space such that $K\df B-A$ is a relatively bounded perturbation of $A$. 

\begin{thm}
\label{Dfy+imS}
Let $f$ be a relatively operator Lipschitz function on $\R$. Then 
the function 
$$
\dg_\rI f=(\dg f)(\y+\ri)
$$
is a Schur multiplier with respect to arbitrary Borel spectral measures;
\end{thm}

%To prove (a), we need the following lemma.

\Pf By Theorem \ref{x+i}, it suffices to show that $(\bx+\ri)f\in\OL(\R)$. Let $A$ and $B$ be bounded self-adjoint operators. We have
\begin{align*}
((\bx+\ri)f)(A)-((\bx+\ri)f)(B)&=f(A)(A+\ri I)-f(B)(B+\ri I)\\[.2cm]
&=(f(A)-f(B))(A+\ri I)+f(B)(A-B).
\end{align*}
Clearly, 
$$
\|f(B)(A-B)\|\le\|f\|_{L^\be}\|A-B\|.
$$
Therefore, $(\bx+\ri)f\in\OL(\R)$ if and only if the following inequality holds:
$$
\|(f(A)-f(B))(A+\ri I)\|\le\const\|A-B\|.
$$
Let $x$ be a vector in our Hilbert space and let $y=(A+\ri I)x$. We have
\begin{align*}
\|(f(A)-f(B))(A+\ri I)x\|&=\|(f(A)-f(B))y\|\le\const\|(A-B)(A+\ri I)^{-1}\|\cdot\|y\|\\[.2cm]
&=\const\|(A-B)(A+\ri I)^{-1}\|\cdot\|(A+\ri I)x\|,
\end{align*}
and so, the function $(\bx+\ri)f$ is operator Lipschitz. $\bl$

Let us now observe that condition \rf{neobkhodimo} does not imply that $f$ is relatively operator Lipschitz.
Indeed, it follows from Theorem \ref{Dfy+imS} and Theorem \ref{x+i} that if $f\in{\rm ROL}$, then 
$(\bx+\ri)f\in\OL(\R)$, and so $f$ must satisfy the following necessary condition (see \cite{APol}):

\begin{thm}
If $f\in{\rm ROL}$, then the function $(\bx+\ri)f\in\OL(\R)$ must belong to the Besov space $B_1^1(\R)$ locally.
\end{thm}

We refer the reader to \cite{APol} for other necessary conditions for a function to be operator Lipschitz.

%{\bf Proof of Statement (b) of Theorem \ref{Dfy+imS}.} Without loss of generality we may assume that $s=0$.
%We have
%$$
%f\big(A_t\big)-f\big(A\big)=
%t\iint_{\R\times\R}\frac{f(x)-f(y)}{x-y}(y+\ri)\,{\rm d}E_{A_t}(x)C\,{\rm d}E_{A}(y), 
%$$
%where the operator $C$ is defined by \rf{C}. If $C=CE_A[-M,M]$, then 
%$$
%f\big(A_t\big)-f\big(A\big)=
%t\iint_{\R\times[-M,M]}\frac{f(x)-f(y)}{x-y}(y+\ri)\,{\rm d}E_{A_t}(x)C\,{\rm d}E_{A}(y)
%$$
%and the same reasoning as in the proof of Lemma \ref{ravnpoM} shows that \rf{proizvodnaya} holds in the strong operator topology. 
%
%To complete the proof, we consider the operators $C_j\df CE_A[-j,j]$ and observe that $C_j\to C$ in the norm.
%$\bl$ 

%\
%
%\section{\bf Relative operator Lipschitzness versus relative trace class Lipschitzness}
%\setcounter{equation}{0}
%
%\

\

\section{\bf Commutator estimates}
\setcounter{equation}{0}
\label{komutator?}

\

In this section we are going to obtain a characterisation of the class of relatively operator Lipschitz functions in terms of commutator and quasi-commutator estimates. 

The following result means, in particular, that a function is relatively operator Lipschitz if and only if it is {\it relatively commutator Lipschitz}, i.e., satisfies statement (b) of the following theorem.

\begin{thm} 
Let $f:\R\to\C$ be a differentiable function. The following
statements are equivalent:

{\em(a)} $f$ is a relatively operator Lipschitz function;

{\em(b)}
$$
\|f(A)R-Rf(A)\|\le\const\|(AR-RA)(A+\ri I)^{-1}\|
$$
for any bounded operator $R$ and any bounded self-adjoint operator $A$;

{\em(c)}
$$
\|f(B)R-Rf(A)\|\le\const\|(BR-RA)(A+\ri I)^{-1}\|
$$
for any bounded operator $R$ and any bounded self-adjoint operators $A$ and $B$;

{\em(d)}
The function $\dg_\rI f$ belongs to $\fM(\R^2)$.

\end{thm}

\Pf First we prove that (a)$\Longrightarrow$(c). 
Let $\f$ be the operator Lipschitz function defined on $\T\setminus\{1\}$ by \rf{na okruzhnosst'} and such that 
$\f(1)=\lim_{|t|\to\be}f(t)$. Then
\bay
\label{UAVB}
\|\f(V)R-R\f(U)\|\le\const\|VR-RU\|
\ey
for every
bounded operator $R$ and for every
unitary operators $U$ and $V$,  see \cite{APol}, Th. 3.1.2.

Let $A$ and $B$  be self-adjoint operators. Put $U=(A-\ri I)(A+\ri I)^{-1}$ and\\ $V=(B-\ri I)(B+\ri I)^{-1}$.
Applying \rf{UAVB} to unitary operators $U$ and $V$ we get
\begin{align*}
\|f(B)R&-Rf(A)\|=\|\f(V)R-R\f(U)\|\le\|\f\|_{\OL(\T)}\|VR-RU\|\\[.2cm]
&=\|\f\|_{\OL(\T)}\|(B+\ri I)^{-1}(B-\ri I)R-R(A-\ri I)(A+\ri I)^{-1}\|\\[.2cm]
&=\|\f\|_{\OL(\T)}\|(B+\ri I)^{-1}\big((B-\ri I)R(A+\ri I)-(B+\ri I)R(A-\ri I)\big)(A+\ri I)^{-1}\|\\[.2cm]
&=2\|\f\|_{\OL(\T)}\|(B+\ri I)^{-1}(BR-RA)(A+\ri I)^{-1}\|\\[.2cm]
&\le2\|\f\|_{\OL(\T)}\|(BR-RA)(A+\ri I)^{-1}\|.
\end{align*}
The implication (c)$\Longrightarrow$(a) is evident.

Let us prove that (b)$\Longrightarrow$(c). Applying (a)
to the self-adjoint operator
$\bs{A}=
\left(
  \begin{array}{cc}
    A &   \0 \\
     \0 &  B \\
  \end{array}
\right)
$
and the bounded operator
$
\bs{R}=\left(
  \begin{array}{cc}
    \0 &  \0  \\
  R  & \0  \\
  \end{array}
\right)$, we find that
$$
\|f(\bs{A})\bs{R}-\bs{R}f(\bs{A})\|\le\const\|({\bf A}\bs{R}-\bs{R}{\bf A})(\bs{A}+\ri\bs{I})^{-1}\|.
$$
It remains to observe that 
$$
\|f(\bs{A})\bs{R}-\bs{R}f(\bs{A})\|=\|f(B)R-Rf(A)\|,
$$
where $\bs{I}\df\left(
  \begin{array}{cc}
    I &   \0 \\
     \0 & I \\
  \end{array}
\right)$,
and
$$
\|(\bs{A}\bs{R}-\bs{R}\bs{A})(\bs{A}+\ri\bs{I})^{-1}\|=\|(BR-RA)(A+\ri I)^{-1}\|.
$$
The implication (c)$\Longrightarrow$(b) is evident.

It remains to observe that the implication (a)$\Longrightarrow$(d) follows from Theorem \ref{otopLi}
and the implication (d)$\Longrightarrow$(a) follows from Theorem \ref{Dfy+imS}. $\bl$

%\begin{thm} 
%Let $f:\R\to\C$ be a differentiable function. The following
%statements are equivalent:
%
%{\em(a)} $f$ is a relatively operator Lipschitz function;
%
%{\em(b)}
%$$
%\|f(A)R-Rf(A)\|\le\const\|(AR-RA)(A+\ri I)^{-1}\|
%$$
%for any bounded operator $R$ and any self-adjoint operator $A$;
%
%{\em(c)}
%$$
%\|f(A)R-Rf(B)\|\le\const\|(AR-RB)(B+\ri I)^{-1}\|
%$$
%for any bounded operator $R$ and any self-adjoint operators $A$ and $B$;
%
%{\em(d)}
%$$
%\|f(A)-f(B)\|\le\const\|(A-B)(A+\ri I)^{-1}\|
%$$
%for any self-adjoint operators $A$ and $B$;
%
%{\em(e)}
%The function $\dg_\rI f$ belongs to $\fM(\R^2)$.
%
%\end{thm}
%
%
%\Pf Let us prove that (a)$\Longrightarrow$(b). Applying (a)
%to the self-adjoint operator
%$\bs{A}=
%\left(
%  \begin{array}{cc}
%    A &   \0 \\
%     \0 &  B \\
%  \end{array}
%\right)
%$
%and the bounded operator
%$
%\bs{R}=\left(
%  \begin{array}{cc}
%    \0 &   R \\
%    \0 &  \0 \\
%  \end{array}
%\right)$, we find that
%$$
%\|f(\bs{A})\bs{R}-\bs{R}f(\bs{A})\|\le\const\|({\bf A}\bs{R}-\bs{R}{\bf A})(\bs{A}+\ri\bs{I})^{-1}\|.
%$$
%It remains to observe that 
%$$
%\|f(\bs{A})\bs{R}-\bs{R}f(\bs{A})\|=\|f(A)R-Rf(B)\|,
%$$
%where $\bs{I}\df\left(
%  \begin{array}{cc}
%    I &   \0 \\
%     \0 & I \\
%  \end{array}
%\right)$
%$$
%\|f(\bs{A})\bs{R}-\bs{R}f(\bs{A})\|=\|(f(A)R-Rf(B))(B+\ri I)^{-1}\|.
%$$
%The implication (b)$\Longrightarrow$(c) is evident.
%
%The implications (c)$\Longrightarrow$(d),  (d)$\Longrightarrow$(a) and (a)$\Longleftrightarrow$(e)
%were proved above. $\bl$

\

\section{\bf Differentiating in the strong operator topology}
\setcounter{equation}{0}
\label{differentsiruem}

\

In this section we prove the differentiability of the function $t\mapsto f(A_t)-f(A)$ in the strong operator topology and obtain a formula for the derivative. This formula will be used in \S\:\ref{fss} to study spectral shift functions in the case of relatively trace class perturbations.

\begin{thm}
\label{sil'no}
Let $f$ be a relatively operator Lipschitz function.
Suppose that a self-adjoint operator $K$ is a relatively bounded perturbation of a self-adjoint operator $A$. Then the function
$t\mapsto f(A_t)-f(A)$ is differentiable at 0 in the strong operator topology and
\begin{align}
\label{proizvodnaya}
\frac{\rm d}{{\rm d}t}\big(f(A_t)-f(A)\big)\Big|_{t=0}&=
\iint_{\R\times\R}\frac{f(x)-f(y)}{x-y}(y+\ri)\,{\rm d}E_{A}(x)C\,{\rm d}E_{A}(y)\nonumber\\[.2cm]
&=\iint_{\R\times\R}\frac{f(x)-f(y)}{x-y}(y^2+1)^{1/2}\,{\rm d}E_{A}(x)G\,{\rm d}E_{A}(y)
\end{align}
where the operator $C$ and $G$ are defined by {\em\rf{C}} and {\em\rf{G}}.
\end{thm}

To prove Theorem \ref{sil'no}, we need the following result.

\begin{thm}
\label{52}
\label{ivsyotakinetreryvno}
Let $A$ be a self-adjoint operator and let $K$ be a self-adjoint relatively bounded perturbation of $A$. Suppose that $\f$ is a continuous function on $\R$ such that the limit $\lim\limits_{|t|\to\be}\f(t)$ exists.
Then the function
$$
\label{AsK}
s\mapsto\f(A+sK)
$$
is continuous in the operator norm.
\end{thm}

\Pf Without loss of generality we may assume that $\lim\limits_{|t|\to\be}\f(t)=0$.
Let $\e>0$. Clearly, there exists an infinitely smooth function $g$ on $\R$ with compact support such that 
$\max\limits_{t\in\R}|\f(t)-g(t)|<\e$. Clearly, $g$ is a relatively operator Lipschitz function and it follows from 
Theorem \ref{otopLi} that
$$
\|g(A+sK)-g(A+s_0K)\|\le\const|s-s_0|\cdot\|K(A+\ri I)^{-1}\|.
$$
The result follows from the inequalities
\begin{multline*}
\|\f(A+sK)-\f(A+s_0K)\|\le\|\f(A+sK)-g(A+sK)\|\\[.2cm]
+\|g(A+sK)-g(A+s_0K)\|+\|g(A+s_0K)-\f(A+s_0K)\|\\[.2cm]
\le2\e+\const|s-s_0|\cdot\|K(A+\ri I)^{-1}\|.\quad\bl
\end{multline*}

To prove Theorem \ref{sil'no}, we need some preparation.

Let $\widehat\R\df\R\cup\{\be\}$ denote the one-point
compactification of the real line $\R$. We recall that each function $f\in\OL(\R$) is
everywhere differentiable on $\widehat\R$ (see Theorem 3.3.3 in \cite{APol}).
Recall that $f$ {\it is said to be differentiable at} $\be$ if there exists a finite limit 
$\lim\limits_{|x|\to\be}x^{-1}f(x)\df f'(\be)$.

\begin{lem}
\label{fipsi}
Let $f$ satisfy the equivalent conditions of Theorem {\em\ref{x+i}}.
Then there are sequences $\{\f_n\}_{n\ge0}$ and $\{\psi_n\}_{n\ge0}$
of continuous functions on $\widehat\R$ and a positive number $M$ such that the following statements hold:

{\em(a)}
$\sum\limits_{n\ge0}|\f_n|^2\le M$ everywhere on $\widehat\R$;

{\em(b)}
$\sum\limits_{n\ge0}|\psi_n|^2\le M$ everywhere on $\widehat\R$;

{\em(c)}
$(x+\ri)(\dg f)(x,y) =\sum\limits_{n\ge0}\f_n(x)\psi_n(y)$
 for all $x,y\in\R$.
 
 Moreover, the constants in {\em(a)} and {\em(b)} depend only on $f$.
\end{lem}

\Pf Note that there exists a finite limit $\lim\limits_{|y|\to\be}f(y)$.
It suffices to observe that $\lim\limits_{|y|\to\be}f(y)=g'(\be)$ for
the operator Lipschitz function $g=(\bx+\ri)f$.

Applying Lemma 3.5.7 in \cite{APol} to the operator Lipschitz function $(\bx+\ri)f$,
we find that
$$
\frac{(x+\ri)f(x)-(y+\ri)f(y)}{x-y}=\sum\limits_{n\ge1}\f_n(x)\psi_n(y)
$$
for some sequences $\{\f_n\}_{n\ge1}$ and $\{\psi_n\}_{n\ge1}$
of continuous functions on $\widehat\R$ such that
$$
\sum\limits_{n\ge1}|\f_n|^2\le\|(\bx+\ri)f\|_{\OL(\R)}\quad\mbox{and}\quad 
\sum\limits_{n\ge1}|\psi_n|^2\le\|(\bx+\ri)f\|_{\OL(\R)}
$$ 
everywhere on $\widehat\R$.

Applying equality \rf{ix}, we otain
$$
 (x+\ri)(\dg f)(x,y) =\sum\limits_{n\ge0}\f_n(x)\psi_n(y)-f(y).
$$
To complete the proof, we put $\f_0(x)=1$, $x\in\R$, and $\psi_0(y)=-f(y)$, $y\in\R$. $\bl$

\medskip

{\bf Proof of Theorem \ref{sil'no}.} The proof is similar to the proof of Theorem 3.5.6 in  \cite{APol}.
We have
\begin{align*}
\frac1t(f(A_t)-f(A))&=
\iint_{\R\times\R}\frac{f(x)-f(y)}{x-y}(y+\ri)\,{\rm d}E_{A_t}(x)K(A+\ri I)^{-1}\,{\rm d}E_A(y)\\[.2cm]
&=\sum_{n\ge0}\f_n(A_t)K(A+\ri I)^{-1}\psi_n(A),
\end{align*}
where $\f_n$ and $\psi_n$  denote the same as in Lemma \ref{fipsi}.
It remains to prove that
$$
\lim_{t\to0}\sum_{n\ge0}\f_n(A_t)K(A+\ri I)^{-1}\psi_n(A)=\sum_{n\ge0}\f_n(A)K(A+\ri I)^{-1}\psi_n(A)
$$
in the strong operator topology. Thus, we need to prove that for an arbitrary vector $u$,
$$
\lim_{t\to0}\sum_{n\ge0}(\f_n(A_t)-\f_n(A))K(A+\ri I)^{-1}\psi_n(A)u=\bf 0
$$
where the series is understood in the sense of  the weak topology of $\mathscr H$ while the limit is taken
in the norm of $\mathscr H$. Assume that $\|u\|=1$.
Let $u_n\df K(A+\ri I)^{-1}\psi_n(A)u$. We have
\begin{align*}
\sum_{n\ge0}\|u_n\|^2&\le\|K(A+\ri I)^{-1}\|^2\sum_{n\ge0}\|(\psi_n(A)u\|^2\\[.2cm]
&=\|K(A+\ri I)^{-1}\|^2\sum_{n\ge0}(|\psi_n|^2(A)u,u)
\le M\|K(A+\ri I)^{-1}\|^2,
\end{align*}
where $M$ is the number in the statement of Lemma \ref{fipsi}.
Let $\e>0$ and choose a natural number $N$ such that $\sum_{n>N}\|u_n\|^2<\e^2$.
 Then it
follows from Lemma 3.5.9 in \cite{APol} that
\begin{align}
\label{khvost}
\Big\|\sum_{n>N}(\f_n(A+tK)-\f_n(A))u_n \Big\|&\le\Big\|\sum_{n>N}\f_n(A+tK)u_n\Big\|\nonumber\\[.2cm]
&+\Big\|\sum_{n>N}\f_n(A)u_n\Big\|
<2M^{1/2}\e
\end{align}
for all $t\in\R$. By Theorem \ref{52},
$$
%\left\|\sum_{n=0}^N(\f_n(A+tK)-\f_n(A))u_n \right\|\le\|K(A+\ri I)^{-1}\|\sum_{n=0}^N\|\f_n(A+tK)-\f_n(A)\|<\e
\left\|\sum_{n=0}^N(\f_n(A+tK)-\f_n(A))u_n \right\|<\e
$$
for all $t$ sufficiently close to zero. Thus, in view of \rf{khvost},
$$
\Big\|\sum_{n\ge0}(\f_n(A+tK)-\f_n(A))u_n\Big\|<(2M^{1/2}+1)\e
$$
for all $t$ sufficiently close to zero. 

The proof of the fact that the derivative on the left of \rf{proizvodnaya} is equal to the second
double operator integral in \rf{proizvodnaya} is the same.
$\bl$

\begin{thm}
Let $f$ be a relatively operator Lipschitz function.
Suppose that a self-adjoint operator $K$ is strictly dominated by a self-adjoint operator $A$. Then the function
$t\mapsto f(A_t)-f(A)$ is differentiable on $\R$ in the strong operator topology and
\begin{align}
\label{proizvodnayas=t}
\!\!\!\!\frac{\rm d}{{\rm d}s}\big(f(A_s)-f(A_t)\big)\Big|_{s=t}\!&=\!
\iint_{\R\times\R}\!\frac{f(x)-f(y)}{x-y}(y+\ri)\,{\rm d}E_{A_t}(x)K(A_t+\ri I)^{-1}\,{\rm d}E_{A_t}(y)\nonumber
\\[.2cm]
&=\iint_{\R\times\R}\!\frac{f(x)-f(y)}{x-y}(y^2+1)^{1/2}\,{\rm d}E_{A_t}(x)K(A_t^2+I)^{-1/2}\,{\rm d}E_{A_t}(y)
\end{align}
for any $t\in(0,1]$.
\end{thm}

\Pf The result follows from Theorem \ref{sil'no} (applied to the operators $K$ and $A_t$) and from Lemma \ref{tret'ya}. $\bl$

\medskip 

In particular, the conclusion of the theorem holds if $K$ is a relatively compact perturbation of $A$.

\

\section{\bf Relatively trace class perturbations and the trace formula}
\setcounter{equation}{0}
\label{fss}

\

In this section we obtain an analogue of the Lifshits--Krein trace formula for relatively trace class perturbation. 
We also describe the maximal class of functions, for which the trace formula is applicable. The main results of this section are based on Theorem \ref{sil'no}.

Let $A$ be a self-adjoint operator and let $K$ be a relatively trace class perturbation of $A$. Suppose that
$f$ is a relatively operator Lipschitz function on $\R$.

The purpose of this section is to prove that the trace formula in the case of relatively trace class perturbations holds for arbitrary relatively operator Lipschitz functions and the corresponding spectral shift function $\bs{\xi}$
satisfies inequality \rf{svesom}. 

\begin{thm}
\label{spectral'no sdvinem}
Let $A$ be a self-adjoint operator and let $K$ be a relatively trace class perturbation of $A$. 
Then there exists a unique measurable real-valued function $\bs\xi$ on $\R$ satisfying
\bay
\label{svesom}
\int_\R\frac{|\bs{\xi}(t)|}{1+|t|}\,{\rm d}t<\be
\ey
and such that for an arbitrary relatively operator Lipschitz function $f$ on $\R$
the following trace formula holds
\bay
\label{nasledil}
\trace\big(f(A+K)-f(A)\big)=\int_\R f'(t)\bs{\xi}(t)\,{\rm d}t.
\ey
\end{thm}

The function $\bs\xi$ satisfying \rf{nasledil} and \rf{svesom} is called the {\it spectral shift function for the pair
$\{A,A+K\}$}.

Note that inequality \rf{svesom} was found earlier in \cite{CS} by different methods. However, in \cite{CS}
trace formula \rf{nasledil} was obtained under considerably more restrictive assumptions on $f$.

Moreover, we describe in this section the maximal class of function $f$, for which the trace formula holds in the case of relatively trace class perturbations, see Theorem \ref{maksimalka} below, and show that this maximal class coincides with the class of relatively operator Lipschitz functions.

The proof of Theorem \ref{spectral'no sdvinem} is close in spirit to the proof of Theorem 6.1 of \cite{PeKr}
and to the proof of Theorem 4.1 of \cite{APun}.

\medskip

\Pf
Consider the one-parametric family of operators $A_t\df A+tK$, $0\le t\le1$. By Theorem \ref{sil'no},
\bay
\label{AQt}
\!\!\frac{\rm d}{{\rm d}s}\big(f(A_s)-f(A_t)\big)\Big|_{s=t}=
\iint_{\R\times\R}\!\!\frac{f(x)-f(y)}{x-y}(y+\ri)\,{\rm d}E_{A_t}(x)K(A_t+\ri I)^{-1}\,{\rm d}E_{A_t}(y),
\ey
where the derivative on the left is understood in the strong operator topology.

We need the following lemma:

\begin{lem}
\label{ogranicheny}
The function $t\mapsto K(A_t+\ri I)^{-1}$ is continuous on $[0,1]$ in the norm of $\bS_1$, and so
$\|K(A_t+\ri I)^{-1}\|_{\bS_1}\le\const$, $0\le t\le1$. 
\end{lem}

\Pf We have
$$
K(A_t+\ri I)^{-1}=K(A+\ri I)^{-1}(A+\ri I)(A_t+\ri I)^{-1}.
$$
Clearly, it suffices to show that the function 
$$
t\mapsto\big((A+\ri I)(A_t+\ri I)^{-1}\big)^{-1}=(A_t+\ri I)(A+\ri I)^{-1}
$$ 
is continuous. This follows immediately from the following obvious equality
$$
(A_t+\ri I)(A+\ri I)^{-1}=I+tK(A+\ri I)^{-1}.\quad\bl
$$

\medskip

Let us complete the proof of Theorem \ref{spectral'no sdvinem}.
Let $Q_t$ be the operator on the right-hand side of \rf{AQt}.
Since the function  $\dg_\rI f$ is a Schur multiplier for arbitrary spectral measures, it follows that
$Q_t\in\bS_1$ for every $t$ in $[0,1]$ and $\sup_t\|Q_t\|_{\bS_1}<\be$.

Since the function $t\mapsto Q_t$ is continuous in the norm of $\bS_1$ by Lemma \ref{ogranicheny},
we obtain
$$
f(A+K)-f(A)=\int_0^1Q_t\,{\rm d}t
$$
in the sense of integration of continuous functions.
Moreover, 
$$
\trace\big(f(A+K)-f(A)\big)=\int_0^1\trace Q_t\,{\rm d}t
$$

We are going to use the following known result (see \cite{PeKr}, Theorem 5.1):

{\it Let $E$ be a Borel spectral measure on a locally compact topological space $X$. Suppose that $\Phi$ is a a Schur multiplier with respect to $E$. If $\Phi$ is continuous in each variable, then
\bay
\label{obflasled}
\trace\left(\iint\Phi(x,y)\,{\rm d}E(x)T\,{\rm d}E(y)\right)=\int\Phi(x,x)\,{\rm d}\mu(x)
\ey
for an arbitrary trace class operator $T$} where the measure $\mu$ is defined by
$$
\mu(\D)\df\trace(E(\D)T).
$$

Applying \rf{obflasled} to the Schur multiplier $\dg_{\rI}f$ and the spectral measure $E_{A_t}$, we find that 
\bay
\label{sledQt}
\trace Q_t=\int_\R f'(x)(x+\ri)\,{\rm d}\nu_t(x),
\ey
where $\nu_t$ is the complex Borel measure on $\R$ defined by 
$$
\nu_t(\D)=\trace\big(E_{A_t}(\D)K(A_t+\ri I)^{-1}\big).
$$

As usual, we identify here the space of regular complex Borel measures on $\R$ with the dual space to the Banach space of continuous functions on $\R$ with zero limit at infinity. Then the function $t\mapsto\nu_t$ is continuous in the weak-$*$ topology on the space of complex Borel measures. 
Indeed, if $h$ belongs to the space
$$
{\rm C}_0(\R)\df\{h:~h~\mbox{is a continuous function on}~\;\R\:~\mbox{and}\:~\lim_{|x|\to\be}h(x) = 0\},
$$
then it is easy to see that
$$
\int h\,{\rm d}\nu_t =\trace(h(A_t)K(A_t+\ri I)^{-1}).
$$
By Theorem \ref{ivsyotakinetreryvno},
the function $t\mapsto h(A_t)$ is a continuous function on [0, 1] in the operator norm, and so
the function $t \mapsto\trace (h(A_t)K(A_t+\ri I)^{-1})$ is continuous.

Consider the real Borel measure $\nu$ on $\R$ defined by
$$
\nu=\int_0^1\nu_t\,{\rm d}t.
$$

We need the following lemma.

\begin{lem}
\label{perkla}
Let $g$ be a bounded function on $\R$ of the first Baire class and let $\mu_t$, $0\le t\le1$, be a parametric family of regular Borel complex measures on $\R$ that depends continuously on the parameter $t$ in the weak-* topology. Then
\bay
\label{gign}
\int_0^1\left(\int_\R g(x)\,\rd\mu_t(x)\right)\rd t=\int_\R g(x)\,\rd\mu(x),
\ey
where $\mu=\int_0^1\mu_t\,\rd t$.
\end{lem}

{\bf The proof of the lemma.} If $g$ is a continuous function and $\lim_{|x|\to\be}=0$, the result is trivial.
In the general case $g$ being on the first Baire class can be approximated by a sequence $\{g_n\}$ of functions in ${\rm C}_0(\R)$ such that 
$$
\sup_{n,x}|g_n(x)|<\be\quad\mbox{and}\quad\lim_ng_n(x)=g(x)\quad\mbox{for all}\quad x\in\R.
$$
Then \rf{gign} holds with $g$ replaced with $g_n$. It remains to pass to the limit as $n\to\be$. $\bl$

\medskip

Let us return to the proof of Theorem \ref{spectral'no sdvinem}. Since relatively operator Lipschitz functions are differentiable everywhere, function $x\to f'(x)(x+\ri)$ is of the first Baire class. Next, this function is bounded on $\R$ because it coincides with the restriction of $\dg_\rI f$ to the diagonal.

Then by \rf{sledQt} and by Lemma \ref{perkla},
$$
\trace\big(f(A+K)-f(A)\big)=\int_\R f'(x)(x+\ri)\,{\rm d}\nu(x).
$$

We have mentioned in \S\:\ref{In} that the condition that $K$ is a relatively trace class perturbation of $A$ implies that
$$
(A+K+\ri I)^{-1}-(A+\ri I)^{-1}\in\bS_1. 
$$

Then the pair $\{A,A+K\}$ has a  spectral shift function $\bs\eta$ such that 
$$
\int_\R\frac{|\bs{\eta}(x)|}{1+x^2}\,{\rm d}x<\be
$$
and
$$
\trace\big(f(A+K)-f(A)\big)=\int_\R f'(x)\bs{\eta}(x)\,{\rm d}x
$$
whenever $f$ is a resolvent Lipschitz function on $\R$. It follows that 
$$
\int_\R f'(x)\bs{\eta}(x)\,{\rm d}x=\int_\R f'(x)(x+\ri)\,{\rm d}\nu(x)
$$
for an arbitrary infinitely smooth function $f$ with compact support. It follows that
$$
\bs{\eta}(x)\,{\rm d}x-(x+\ri)\,{\rm d}\nu(x)=c\,{\rm d}x
$$
for some constant $c$. In particular, this means that $\nu$ is absolutely continuous with respect to Lebesgue measure. We can define now the function $\bs{\xi}$ by
$$
\bs{\xi}(s)\df(s+\ri)\frac{{\rm d}\nu(x)}{{\rm d}x}(s).
$$
Clearly, $\bs\xi$ satisfies \rf{svesom} and equality \rf{nasledil} holds.

Suppose now that $f$ is a real-valued relatively operator Lipschitz function. Clearly, in this case
$$
\trace\big(f(A+K)-f(A)\big)\in\R.
$$
It follows that
$$
\trace\big(f(A+K)-f(A)\big)=\int_\R f'(x)\bs{\xi}(x)\,{\rm d}x=\int_\R f'(x)(\re\bs{\xi})(x)\,{\rm d}x.
$$
Thus, both $\bs\xi$ and $\re\bs{\xi}$ are spectral shift functions for the pair $\{A,A+K\}$ and both of them satisfy inequality \rf{svesom}. To show that $\bs\xi$ is real-valued, it suffices to show that a spectral shift function satisfying \rf{svesom} is unique.

Indeed, if $\bs{\xi}_1$ and $\bs{\xi}_2$ are spectral shift functions such that
$$
\int_\R\frac{|\bs{\xi}_1(x)|}{1+|x|}\,{\rm d}x<\be\quad\mbox{and}\quad
\int_\R\frac{|\bs{\xi}_2(x)|}{1+|x|}\,{\rm d}x<\be,
$$
then $\bs{\xi}_1-\bs{\xi}_2=c$ is a constant function. Then
$$
\int_\R\frac{|\bs{\xi}_1(x)|}{1+x^2}\,{\rm d}x<\be\quad\mbox{and}\quad
\int_\R\frac{|\bs{\xi}_2(x)|}{1+x^2}\,{\rm d}x<\be.
$$
Since $\bs\xi_1$ and $\bs\xi_2$ are spectral shift functions that corresponds to the resolvent trace class perturbation $K$, $\bs{\xi}_1-\bs{\xi}_2=c$ is a constant function. This is a well known fact. It follows  from the fact that an integrable spectral shift function for a pair of unitary operators with trace class difference is unique modulo a constant additive.

If the constant $c$ were nonzero, it would follow that
$$
\int_\R\frac{|c|}{1+|x|}\,{\rm d}x<\be.\quad\bl
$$

The following result shows that the class ${\rm ROL}$ is the maximal class of functions on $\R$ for
which trace formula \rf{nasledil} holds.

\begin{thm}
\label{maksimalka}
Let $f$ be a differentiable function on $\R$ such that trace formula {\em\rf{nasledil}}
holds for arbitrary self-adjoint operators $A$ and $K$ such that $K$ is a relatively trace class perturbation of $A$. Then $f\in{\rm ROL}$.
\end{thm}

\Pf Indeed, it follows from Theorem \ref{rolryal} that if $f\not\in{\rm ROL}$, then there are self-adjoint operators $A$ and $K$ such that $K$ is a relatively trace class perturbation of $A$ but $f(A+K)-f(A)\not\in\bS_1$. $\bl$

%\
%
%\
%
%\
%
%
%It follows from +++ that $\nu$ is absolutely continuous with respect to Lebesgue measure. Let 
%${\rm d}\nu(x)=\bs\eta(x){\rm d}x$. Then $\bs\eta\in L^1$. Put $\bs\xi(x)=(x^2+1)^{1/2}\bs\eta(x)$. Clearly,
%$$
%\int_\R\frac{|\bs\xi(x)|}{1+|x|}\,{\rm d}(x)<\be
%$$
%and
%$$
%\trace\big(f(A+K)-f(A)\big)=\int_0^1f'(x)\bs\xi(x)\,{\rm d}(x)
%$$
%for 
%Note that the above reasoning is similar to the one in the proof of Theorem 6.1 of \cite{PeKr}.

\

\

\
 
\begin{footnotesize}
 
\noindent
\begin{tabular}{p{8cm}p{15cm}}
A.B. Aleksandrov & V.V. Peller \\
St.Petersburg State University & St.Petersburg State University \\
Universitetskaya nab., 7/9  & Universitetskaya nab., 7/9\\
199034 St.Petersburg, Russia & 199034 St.Petersburg, Russia \\
\\

St.Petersburg Department &St.Petersburg Department\\
Steklov Institute of Mathematics  &Steklov Institute of Mathematics  \\
Russian Academy of Sciences  & Russian Academy of Sciences \\
Fontanka 27, 191023 St.Petersburg &Fontanka 27, 191023 St.Petersburg\\
Russia&Russia\\
email: alex@pdmi.ras.ru& email: peller@math.msu.edu
\end{tabular}
\end{footnotesize}

\end{document}

%% file: classstartscr.tex
\newcommand{\lb}{\linebreak}

\newcommand{\e}{\varepsilon}

\newcommand{\z}{\zeta}

\newcommand{\f}{\varphi}

\newcommand{\D}{\Delta}
\renewcommand{\L}{\Lambda}

\newcommand{\E}{{\mathscr E}}

\newcommand{\C}{{\Bbb C}}
\newcommand{\T}{{\Bbb T}}

\newcommand{\R}{{\Bbb R}}

\newcommand{\0}{{\boldsymbol{0}}}

\newcommand{\bs}{\boldsymbol}

\newcommand{\bS}{{\boldsymbol S}}

\newcommand{\rf}[1]{(\ref{#1})}

\newcommand{\df}{\stackrel{\mathrm{def}}{=}}

\newcommand{\re}{\operatorname{Re}}

\newcommand{\trace}{\operatorname{trace}}
\newcommand{\rank}{\operatorname{rank}}
\newcommand{\const}{\operatorname{const}}

\newcommand{\eeq}{\end{equation}}
\newcommand{\beq}{\begin{equation}}
\newcommand{\bay}{\begin{eqnarray}}
\newcommand{\ba}{\begin{align*}}
\newcommand{\ea}{\end{align*}}
\newcommand{\ey}{\end{eqnarray}}
\newcommand{\bey}{\begin{eqnarray*}}
\newcommand{\eey}{\end{eqnarray*}}

\newcommand{\be}{\infty}

\newcommand{\bl}{\blacksquare}

\newcommand{\Range}{\operatorname{Range}}

\newcommand{\Pf}{{\bf Proof. }}

\renewcommand{\re}{\operatorname{Re}}

\newtheorem{thm}{\hspace{\parindent}Theorem}[section]

\newtheorem{cor}[thm]{\hspace{\parindent}Corollary}
\newtheorem{lem}[thm]{\hspace{\parindent}Lemma}

%% file: Vsyo_otnositel_no_.bbl
\begin{thebibliography}{99}
\label{bibl}

%\bibitem[AP1]{AP2}  {\sc A.B. Aleksandrov} and {\sc V.V. Peller},  {\it Operator H\"older--Zygmund functions}, Advances in Math.
%{\bf224} (2010), 910--966.

%\bibitem[AP3]{AP3}  {\sc A.B. Aleksandrov} and {\sc V.V. Peller},  {\em Functions of operators under perturbations of
%class $\bS_p$}, J. Funct. Anal. {\bf258} (2010), 3675--3724.

\bibitem[A]{A} A. B. Aleksandrov, Operator Lipschitz functions and linear-fractional transformations,
Zapiski Nauchn. Semin. POMI {\bf401} (2012), 5--52 (Russian).
English transl.: J. Math. Sci. (N. Y.) {\bf194:6} (2013), 603--627.

\bibitem[AP1]{APol}  {\sc A.B. Aleksandrov} and {\sc V.V. Peller}, {\it Operator Lipschitz functions}, Uspekhi Mat. Nauk
{\bf71:4} (2016), 3--106 (Russian).
English transl.: Russian Math. Surveys {\bf71:4}, 605--702.

\bibitem[AP2]{APun} {\sc A.B. Aleksandrov} and {\sc V.V. Peller}, {\it Krein's trace formula for unitary operators and operator Lipschitz functions}, Funkts. Anal. i Ego Pril. {\bf50:3} (2016), 1--11 (Russian).
English transl.: Funct. Anal. Appl. {\bf50:4} (2016), 167--175.

\bibitem[AP3]{AP3} {\sc A.B. Aleksandrov} and {\sc V.V. Peller}, {\it Haagerup tensor products and Schur multipliers}, Algebra i Analiz {\bf36:5} (2024), 

%\bibitem[APPS1]{APPS} {\sc A.B. Aleksandrov, V.V. Peller, D. Potapov}, and
%{\sc F. Sukochev}, {\em Functionsof perturbednormal operators}, C.R. Acad. Sci. Paris, S\'er I {\bf348} (2010), 553--558.
%
%\bibitem[APPS2]{APPS2} {\sc A.B. Aleksandrov, V.V. Peller, D. Potapov}, and
%{\sc F. Sukochev}, {\em Functions of normal operators under perturbations},
%Advances in Math. {\bf226} (2011), 5216-–5251.

%\bibitem[Be]{Be} {\sc G. Bennett}, {\em Schur multipliers}, Duke Math. J. {\bf44} (1977), 603--639.

\bibitem[BS1]{BS1} {\sc M.S. Birman} and {\sc M.Z. Solomyak},
{\em Double Stieltjes operator integrals},
Problems of Math. Phys., Leningrad. Univ. {\bf1} (1966), 33--67 (Russian).
English transl.: Topics Math. Physics {\bf1} (1967), 25--54, Consultants Bureau Plenum
Publishing Corporation, New York.

\bibitem[BS2]{BS2} {\sc M.S. Birman} and {\sc M.Z. Solomyak},
 {\em Double Stieltjes operator integrals. II},
 Problems of Math. Phys., Leningrad. Univ. {\bf2} (1967), 26--60 (Russian).
English transl.: Topics Math. Physics {\bf2} (1968), 19--46, Consultants Bureau Plenum
Publishing Corporation, New York.

\bibitem[BS3]{BSsl} {\sc M.S. Birman} and {\sc M.Z. Solomyak}, {\em Remarks on the spectral
shift function},  Zapiski Nauchn. Semin. LOMI {\bf27} (1972), 33--46 (Russian).

English transl.: J. Soviet Math. {\bf3} (1975), 408--419.

\bibitem[BS4]{BS4} {\sc M.S. Birman} and {\sc M.Z. Solomyak},
{\em Double Stieltjes operator integrals. III},
Problems of Math. Phys., Leningrad. Univ. {\bf6} (1973), 27--53 (Russian).

\bibitem[BS5]{BSSst} {\sc M.Sh. Birman and M.Z. Solomyak}, {\it Spectral theory of
selfadjoint operators in Hilbert space},
Mathematics and its Applications (Soviet Series),
D. Reidel Publishing Co., Dordrecht, 1987.

\bibitem[BS6]{BStp} {\sc M.S. Birman} and {\sc M.Z. Solomyak},
{\em Tensor product of a finite number of spectral measures is always a spectral measure}, Integral Equations Operator Theory {\bf24} (1996),  179--187.

\bibitem[CS]{CS} {\sc A. Chattopadhyay} and {\sc A. Skripka}, {\it Trace formulas for relative Schatten class perturbations}, J. Funct. Anal. {\bf274} (2018), 3377--3410.



\bibitem[DK]{DK} {\sc Yu.L. Daletskii} and {\sc S.G. Krein}, {\em Integration and differentiation of
functions of Hermitian operators and application to the theory of perturbations} (Russian), Trudy Sem.
Functsion. Anal., Voronezh. Gos. Univ. {\bf1} (1956), 81--105.



\bibitem[F]{F2}  {\sc Yu.B. Farforovskaya}, {\em An example of a Lipschitzian function of selfadjoint
operators that yields a nonnuclear increase under a nuclear perturbation}.  Zap. Nauchn. Sem.
Leningrad. Otdel. Mat. Inst. Steklov. (LOMI)  {\bf30}  (1972), 146--153 (Russian).

\bibitem[GK]{GK}
I.C. Gokhberg and M.G. Krein, {Introduction to the theory of linear nonselfadjoint
operators in Hilbert space}, "Nauka", Moscow, 1965;\\  English transl., Amer. Math. Soc.,
Providence, R.I., 1969.



%\bibitem{JW} {\sc B.E. Johnson} and {\sc J.P. Williams}, {\em The range of a normal derivation},
%Pacific J. Math. {\bf58} (1975), 105--122.


\bibitem[Kr]{Kr} {\sc M.G. Krein}, {\em On a trace formula in perturbation theory},
Mat. Sbornik {\bf33} (1953), 597--626 (Russian).

%\bibitem[Kr2]{Kr2} {\sc M.G. Krein}, {\em Some new studies in the theory of perturbations of self-adjoint operators}, First Math. Summer School, Part I (Russian), 103-–187, Naukova Dumka, Kiev, 1964.

\bibitem[L]{L} {\sc I.M. Lifshits}, {\em On a problem in perturbation theory
connected with quantum statistics}, Uspekhi Mat. Nauk {\bf7} (1952), 171--180 (Russian).

\bibitem[MNP]{MNP}
{\sc M.M. Malamud, H. Neidhardt} and {\sc V.V. Peller}, {\it Absolute continuity of spectral shift},
J. Funct. Anal. {\bf276} (2019), 1575--1621.

\bibitem[Pee]{Pee} {\sc J. Peetre},
{\em New thoughts on Besov spaces}, Duke Univ. Press., Durham, NC, 1976.


%\bibitem{Pe1} {\sc V.V.Peller}, {\em Hankel operators of class ${\bf S}_{p}$
%and their applications (rational approximation, Gaussian processes,
%the problem of majorizing operators)}, Mat. Sbornik,
%{\bf 113} (1980), 538-581.
%
%English Transl. in Math. USSR Sbornik, {\bf 41}
%(1982), 443-479.


\bibitem[Pe1]{Pe1} {\sc V.V. Peller},
{\em Hankel operators in the theory of perturbations of unitary and self-adjoint operators},
Funktsional. Anal. i Prilozhen. {\bf19:2}  (1985),
37--51 (Russian).
English transl.: Funct. Anal. Appl. {\bf19} (1985) , 111--123.

%\bibitem[Pe]{Pe0} {\sc V.V. Peller}, {\em For which $f$ does $A-B\in{\bf S}_{p}$
%imply that $f(A)-f(B)\in{\bf S}_{p}$?}, Operator Theory, Birkh\"{a}user,
%{\bf 24} (1987), 289-294.


\bibitem[Pe2]{Pe2} {\sc V.V. Peller},
{\em Hankel operators in the perturbation theory of of unbounded self-adjoint operators}.
Analysis and partial differential equations,  529--544,
Lecture Notes in Pure and Appl. Math., {\bf122}, Dekker, New York, 1990.

%\bibitem{Pe4} {\sc V.V. Peller}, {\em Hankel operators and their applications,}
%Springer-Verlag, New York, 2003.

%\bibitem{Pe5} {\sc V.V. Peller}, {\em Multiple operator integrals in perturbation theory}, Bull. Math. Sci. {\bf6} (2016), 15--88.

\bibitem[Pe3]{PeKr}{\sc V.V. Peller}, {\em The Lifshits--Krein trace formula and operator Lipschitz functions}, Proc. Amer. Math. Soc. {\bf144} (2016), 5207--5215.






\end{thebibliography}
